\documentstyle[12pt]{article}
\input amssym.tex
\normalbaselines
\begin{document}
\newcommand{\EE}{{\cal E}}
\newcommand{\PP}{{\cal P}}
\newcommand{\Pe}{\PP_\eps}
\newcommand{\Po}{\overline{\PP}}
\newcommand{\Peo}{\overline{\Pe}}
\newcommand{\II}{{\cal I}}
\newcommand{\JJ}{{\cal J}}
\newcommand{\OG}{{\cal O}_q(G)}
\newcommand{\LL}{{\cal L}}
\newcommand{\QQ}{{\cal Q}}
\newcommand{\C}{\Bbb C}
\newcommand{\OO}{{\Bbb O}_{st}}
\newcommand{\eps}{\varepsilon}
\newcommand{\Z}{\Bbb Z}
\newcommand{\Q}{\Bbb Q}
\newcommand{\mod}{\mbox{mod}}
\newcommand{\Spec}{\mbox{Spec}}
\newcommand{\Goldie}{\mbox{Goldie}}
\newcommand{\Center}{\mbox{Center}}
\newcommand{\Disc}{\mbox{Disc}}
\newcommand{\Ad}{\mbox{Ad}}
\newcommand{\Ker}{\mbox{Ker}}
\newcommand{\Fract}{\mbox{Fract}}
\newcommand{\N}{\Bbb N}
\newcommand{\M}{\Bbb M}
\newcommand{\xt}{\tilde{x}}
\newcommand{\la}{\lambda}
\newcommand{\pht}{\tilde{\phi}}
\author{A.~N.~Panov.}
\date{\it{Russia,443011,Samara,ul.Akad.Pavlova 1, Samara State University,
Mathematical Department, E-mail: panov@info.ssu.samara.ru}}
\title{
Stratification of prime spectrum of quantum solvable algebras.}
\maketitle
\footnote{The work is supported by RFFI grant  99-01-00014.}
{\bf Introduction.}
We consider the class of Noetherian ring, appeared as
a result of quantization of algebraic groups and their
representations within framework of mathematical physics.
One set up the problem of description
of prime and primitive spectrum of these rings.
This problem has been solved first for algebras of low dimension, 
then for the case
$GL_q(n)$, later for general case of regular functions on
quantum semisimple group [J].
Simultaneously some authors consider the examples
of Quantum Weyl algebra, Quantum ortogonal and symplectic spaces,
Quantum Heisenberg algebra.
One can present new examples considering subalgebras, fartor algebras
of the above algebras and also their deformations.
It was noted that prime ideals in these algebras
have common properties like: 
all prime ideals are comletely prime, there exists a stratification of 
spectrum by locally closed sets,
the hypothesis on  fields of fractions of prime factors is true 
(quantum Gelfand-Kirillov conjecture) and others.

In the survey [G] all known examples are considered,
the common properties are extracted 
and the problem of creating the general theory covering all these examples
is setted up.
The classical analog of this proposed theory
is the theory of 
prime ideals in universal enveloping algebra of a solvable Lie algebra
[D],[BGR].

In this paper the author presents his version of the system of axioms, 
covering the listed examples.  
First, all mentioned algebras are solvable (see Definition 1.2).
The class of solvable algebras obeying Conditons Q1-Q3
covers "quantum" algebras, "classical" algebras
(universal enveloping algebra of solvable Lie algebra)
and there are "mixed" algebras [AD].
 
Considered by the author Condition Q4 may be treated as some condition for
quantum  algebra to be "algebraic".
The classical analog of this property is as follows.
The Lie algebra of an algebraic group, defined over  ${\Bbb Z}$,
admits $p$-structure after reduction modulo prime number $p$.
Its universal enveloping algebra becomes 
finite dimensional over the center.
Condition Q4 helps to isolate "quantum" algebras (Definition 2.8).

All Conditions Q1-Q4 are easily checkable.
The Condition Q4 in the examples is derived from the property
that quantum algebras are finite over center at roots of 1.

In the paper the finite stratification of spectrum of quantum solvable algebra,
obeying Conditions Q1-Q4, is constructed (Theorems 3.4, 3.10)
The hypothesis on skewfields of fractions is proved for prime factors
of these algebras (Corollary 3.8).
Earlier
this hypothesis was proved for the full skew field of fractions of
the algebra [P1]. Using another conditions and methods,
the stratification of quantum solvable algebras was obtained in [C1,C2].

The author is thankfull for J.Alev, G.Cauchon for discussions.

\section{Basic definitions and examples}

Consider a commutative noetherian 
domain $C$ with the property $1+1+\cdots+1\ne 0$.
Let $F$ be the field of fractions of $C$. The characteristic of $F$ is zero.\\
{\bf Definition 1.1} Let $R$ be a ring and $C$ contains in the center of $R$. 
( we shall refer $R$ as $C$-algebra). We say that two elements $a,b\in R$
$q$-commutes, if there holds $ab=qba$ for some invertible element
$q\in C$.\\
Let $Q=(q_{ij})$ be a $n\times n$-matrix.
Suppose that the entires $q_{ij}$ are invertible elements of $C$ and
$q_{ij}q_{ji}=q_{ii}=1$.\\
{\bf Definition 1.2}.
We say that a ring $R$ is quantum solvable over $C$, if
$R$ is generated by the elements 
$x_1, x_2,\ldots, x_n, k_1^{\pm 1},\ldots,
k_m^{\pm 1}$ 
such that the monomials $x_1^{t_1}\cdots x_n^{t_n}k_1^{s_1}\cdots
k_m^{s_m}$ with $t_i\in\N$,$s_i\in\Z$ form a free basis over $C$ with the 
relations\\
1) the elements $k_1,\ldots,k_m$ $q$-commute with all generators  
$x_1, x_2,\ldots, x_n, k_1^{\pm 1},\ldots,
k_m^{\pm 1}$ ;\\
2) for $i<j$ the holds 
$$x_ix_j=q_{ij}x_jx_i+r_{ij},\eqno (1.1)$$
where $r_{ij}$ is the element of subalgebra
$R_{i+1}$ generated by $x_{i+1},\ldots, x_n, k_1^{\pm 1},\ldots,
k_m^{\pm 1}$ .

There exists the chain of subalgebras
$$R=R_1\supset R_2\supset\cdots\supset R_n\supset 
R_{n+1}= C_q[k_1^{\pm 1},\ldots,
k_m^{\pm 1}]\supset C.$$
Here $R_{n+1}$ is an algebra of twisted Laurent polynomials (Definition 3.1).
A quantum solvable algebra is an iterated skew extension of $C$
(see [MC-R],[GL1]).
This means that  
 for all $i\in [1,n]$ the map 
$\tau_i: x_j\mapsto q_{ij}x_j, i>j$ is extended to automorphism of
$R_{i+1}$ and the map $\delta_i:x_j\mapsto r_{ij}$ is extended to 
$\tau_i$-derivation of $R_{i+1}$. We have  
 $\delta_i(ab)=\delta_i(a)b+\tau_i(a)\delta_i(b)$ and
$x_ia=\tau_i(a)x_i+\delta_i(a)$ for all $a,b\in R_{i+1}$.
This extension denoted by $R_i=R_{i+1}[x;\tau_i,\delta_i]$ and is called
skew extension or Ore extension ([MC-R],1.2,1.9.2).
All automorphisms $\tau_i$ are identical on $C$ and all 
$\tau_i$-derivations $\delta_i$ are equal to zero on $C$.
A quantum solvable algebra is a Noetherian domain ([MC-R],1.2.9).
We put some more conditions. \\
{\bf Condition Q1}. 
We claim that $\delta_i\tau_i=q_i\tau_i\delta_i$ for some
invertible elements $q_i\in C$ 
(if $\delta_i=0$, we put $q_i=1$). 
So $R$ is an iterated $q$-skew extension of $C$ ([GL1,2].
\\
{\bf Remark}. The condition Q1 implies that for every pair
$R_i\supset R_{i+1}$ the automorphism
$\tau_i:R_{i+1}\mapsto R_{i+1}$ can be extended to
automorphism of $R_i$ via $\tau_i(x_i)=q_i^{-1}x_i$.\\
{\bf Condition Q2}. 
Denote 
by $\Gamma$ the subgroup generated by all $q_{ij},q_i$
in the group $C^*$ of invertible elements in $C$.
We claim that $\Gamma$ is torsionfree.\\
{\bf Remark}. 
If $R$ is quantum solvable algebra obeying Q1,Q2 and
$\II$ is a prime ideal in $R$ 
such that the image of $\Gamma$ in $C/(I\cap C)$ is torsionfree, then
$\II$ is completely prime 
([GL2],Theorem 2.3). In particular, all prime ideal with zero intersection with
$C$ are completely prime.\\
{\bf Condition Q3}. All automorphisms $\tau_i$ are extended to diagonal
automorphisms of $R$ (i.e the monomials 
$x_1^{t_1}\cdots x_n^{t_n}k_1^{s_1}\cdots
k_m^{s_m}$, $t_i\in\N$,$s_i\in\Z$ are eigenvectors) with 
eigenvalues in $\Gamma$. 
Denote by $H$ the group of diagonal automorphisms 
generated by $\tau_1,\tau_2,\ldots,\tau_n$.
We shall refer eigenvalues as $H$-weights.

  Let $\EE$ be a prime ideal in $C$. We denote $\eps:C\mapsto C/\EE$
and call $\eps$ a specialization
of $C$. We denote $R_\eps=R/(R\EE)$.\\
{\bf Condition Q4}. $R$ is a pure $C$-algebra in the sense of [P1].\\  
It means that the set 
$$\Omega=
\{\EE=Ker(\eps) \vert\quad R_\eps\quad\mbox{is a domain and finite over its 
center}\}$$
is dense in $Spec(C)$ in Zariski topology. 
Notice that in the case $R$ is quantum solvable the algebras 
$R_\eps$ are always domains.\\ 
We shall give some examples of quantum solvable algebras 
obeying Conditions Q1-Q4.
{\bf Example 1. Quantum matrices}.\\
Consider the commutative algebra $C=\C[q_{ij}^{\pm},c^{\pm}: i<j]$.
Let $P,Q\in Mat_n(k)$ be matrices such that
$p_{ij}q_{ij}= c^{sgn(j-i)}$, $p_{ij}p_{ji}=q_{ij}q_{ji}=p_{ii}=q_{ii}=1$.
The algebra $M_{P,Q,c}(n)$ of regular functions on quantum matrices
is generated by $C$ and the entries of quantum matrix  
$\{a_{ti}\}_{t,i=1}^n$. The  elements of $C$ lie in the center.
The quantum entries obey the relations 
$$a_{ti}a_{sj}- q_{ts}q_{ij}^{-1}a_{sj}a_{ti}= 
(p_{ts}^{-1}- q_{ij})a_{si}a_{tj}$$
for $i<j$, $t<s$ or $i>j$, $t>s$ and
$$ a_{ti}a_{sj}= p_{ts}^{-1}q_{ij}^{-1}a_{sj}a_{ti} $$
in all other cases (see [MP]). 
Denote $x_{(t-1)n+i}=a_{ti}$ and $R_m$ is the 
algebra generated by $x_j$, $j=1,\ldots,m$. Consider the filtraton:
$$R=R_{n^2}\supset \ldots \supset R_1\supset R_{0}= C.$$
The $M_{P,Q,c}(n)$ is a quantum solvable algebra ([De],Theorem 2.22) 
with respect to this filtration.
and a pure $C$-algebra (see [P1]).
The group $\Gamma$ is torsionfree.
The subring $R_m$ is the extension $R_m=R_{m-1}[a_{ti};\tau_{ti},\delta_{ti}]$
where $m=(t-1)n+i$, $\tau_{ti}$ is the automorphism of $R_{m-1}$ 
and $\delta_{ti}$ is $\tau_{ti}$-derivation of $R_{m-1}$ with
$\delta_{ti}\tau_{ti}=c^{-1}\tau_{ti}\delta_{ti}$.
Here 
$$
\tau_{ti}(a)=\left\{\begin{array}{l}
p_{st}^{-1}a_{si},\quad {\mbox{if}}\quad a=a_{si}, s<t;\\
q_{ij}^{-1}a_{tj},\quad {\mbox{if}}\quad  a=a_{tj}, j<i;\\
c a_{ti},\quad {\mbox{if}}\quad a=a_{ti} s<t;\\
p_{st}^{-1}q_{ij}^{-1}c^{-1}a,\quad {\mbox{in all other cases}}.\\
\end{array}\right.
$$
The algebra $M_{P,Q,c}(n)$  admits automorfisms via multiplication of columns
and rows of quantum matrix by elements of $C^*$.
One can extend $\tau_{ti}$ to automorphism of $M_{P,Q,c}(n)$
as the product
of the following automorphisms:
multiplication of $s$-row on $p_{ts}^{-1}c^{-\frac{1}{2}}$ for all $s$,
besides $s=t$;
multiplication of $j$-row on $q_{ij}^{-1}c^{-\frac{1}{2}}$ for all $j$,
besides $j=i$;
multiplication of $t$-row and $i$ column on $c^{\frac{1}{2}}$.
Hence  $M_{P,Q,c}(n)$  is quantum solvable algebra obeying Conditions Q1-Q4.
One can find new examples, considering subalgebras and factor algebras of
quantum matrices. The triangular quantum matrices and 
quantum $(m,n)$-matrices are subalgebras
of $M_{P,Q,c}(n)$ are also the examples.
More  generally, consider some subset 
$W$ in the set of pairs $(k,l)$,
$k,l=\overline{1,n}$ such that if two pairs 
$(s,j),(t,i)$ with $s<t,j<i$ lie in $W$,
then
$(s,i),(t,j)$ also lie $W$. One can treat subalgebra
$M_{P,Q,c}(W)$
generated by entries $a_{ti}, (t,i)\in W$
and it is an example.\\
{\bf Example 2. Quantum Weyl algebra}.\\
 The matrices $P,Q$ and $c$ as in Quantum matrices.
A Quantum Weyl algebra $A_{P,Q,c}(n)$ 
is generated by $C$ and $x_i$,$y_i$, $i=\overline{1,n}$ with the relations\\ 
$x_ix_j=q_{ij}^{-1}x_jx_i$;\\ 
$y_iy_j=p_{ij}y_jy_i$ ;\\
$y_ix_j=p_{ji}x_jy_i$ for $i<j$;\\ 
$y_ix_j=q_{ij}x_jy_i$ for $i>j$;\\
$y_ix_i=1+c^{-1}x_iy_i+(c^{-1}-1)\sum_{\alpha=i+1}^nx_\alpha y_\alpha$; 
([De],4.2).

Consider the $C$-basis $y_1,\ldots,y_n,x_n,\ldots,x_1$ 
and the chain of subalgebras $R_i=<y_i,y_{i+1},\ldots,y_n,x_n,\ldots,x_1>$:
$$R=R_1\supset R_2\supset\ldots\supset R_{n+1}= C_q[x_1,\ldots,x_n]$$
Here $R_i=R_{i+1}[y_i,\tau_i,\delta_i]$ with
$\delta_i\tau_i=c^{-1}\tau_i\delta_i$.
The Quantum Weyl algebra 
is an iterated $q$-skew polynomial extension of $C$. 
The Quantum Weyl algebra is quantum solvable and pure $C$-algebra [P1].
The automorphism $\tau_i$ acts on $R_{i+1}$ as follows
$$
\tau_{i}(a)=\left\{\begin{array}{l}
p_{ij}y_j,\quad {\mbox{if}}\quad a=y_j, i<j;\\
p_{ji}x_j,\quad {\mbox{if}}\quad a=x_j, i<j;\\
q_{ij}x_j,\quad {\mbox{if}}\quad a=x_j, i>j;\\
c^{-1}x_i,\quad {\mbox{if}}\quad a=x_i.
\end{array}\right.
$$
The Quantum Weyl algebra admits torus automorphisms 
$x_j,y_j\mapsto a_jx_j,a_j^{-1}y_j$.
Putting $\tau_i(y_j)=q_{ji}$ for $j<i$ and $\tau_i(y_i)=cy_i$,
we extend $\tau_i$ to an automorphism of 
$A_{P,Q,c}(n)$. Conditions Q1-Q4 are true.\\
{\bf Example 3}. $U^+_q=U_q({\bf n}^+)$ is quantum nilpotent
universal enveloping algebras ([G]). 
over $ C=\C[q,q^{-1}, (q^{d_i}-q^{-d_i})^{-1}]$
by $E_i$, $i=\overline{1,n}$ with the quantum 
Chevalley-Serre relations.
Consider the set of simple roots $\pi=\{\alpha_1,\ldots,\alpha_n\}$, 
the set of positive roots $\Delta ^+$. 
Fix a reduced expression $w_0=s_{i_1}\ldots s_{i_N}$ of the longest element
in the Weyl group W. Consider the following convex ordering 
$$\beta_1=\la_{i_1}, \beta_2=s_{i_1}(\la_2),\ldots, 
\beta_N= s_{i_1}\ldots s_{i_{N-1}}(\la_N)$$
in $ \Delta^+$.
Introduce the corresponding $root~vectors$ for $s=\overline{1,N}$ (see[Lu]):
  $$ E_{\beta_s}=T_{i_1}\cdots T_{i_{s-1}}E_{i_s}.$$
For $m=(m_1,\ldots,m_N)$  let 
$E^m= E_{\beta_1}^{m_1}\cdots E_{\beta_N}^{m_N}.$
The elements $E^m$, $m\in \Z_+^N$ form a basis of $U^+_q$
over $C$ [L]. The algebra $U^+_q$ is quantum solvable . Moreover,
it is an iterated $q$-skew extention ([G],2.2) with respect to 
the sequence of extentions $R_i=R_{i+1}[E_{\beta_i};\tau_i,\delta_i]$
(see [G]. 2.2).
Here automorphism $\tau_i$ of $R_{i+1}$
is presented
$\tau_i(E_\gamma)= q^{-(\beta_i,\gamma)}E_\gamma$
and extends to an automorphism $U^+_q$.
The algebra $U^+_q$ is a pure $C$ algebra ([P1]).
 All this proves that $U^+_q$ obey Q1-Q4.

Universal quantum Borel subalgebra
$U_q({\bf b}^+)$ is generated by diagonal torus
$T=\C[K^{\pm 1}_1,\ldots,K^{\pm 1}_n]=$ and $U_q({\bf n}^+)$.
It is also an example.
Considering the multiparameter versions of these algebras and
some subalgebras in them, we get new examples.\\
{\bf Example 4}. The algebra $\OG$ of regular functions
on quantum semisimple Lie group $G$(see [J]).

Let $U_q({\bf g})$ be an quantum universal enveloping algebra of semisimple 
complex Lie algebra ${\bf g}$.   
Let $C=\C[q,q^{-1}]$ and $F=\Fract(C)=\C(q).$
We consider $U_q({\bf g})$ as an algebra over $F$.
By definition, $\OG$ is a Hopf dual algebra  of $U_q({\bf b})$.
 The algebra $\OG$ is generated by the  matrix elements of finite dimensional
representations  of $\OG$.
Consider the extension of diagonal torus $T$ by $\tilde{T}$ generated by the elements
$K_\gamma$,$\gamma\in P(\pi)$ (here $P(\pi)$ is the lattice of weights).  
Consider two extensions of Borel subalgebras
$\tilde{V}_q^\pm=\tilde{T}U_q^\pm$.

There exists the denominator set $S=c_e$ (notation of [J]) such that
every prime ideal in $\OG$ have zero intersection with
$S$ ([J],9.3.10). Studing prime ideals, we can substitute $\OG$ by
its localizaton $\OG_S$. 
On other hand  ([J],9.2.14), the algebra $\OG_S$ is  
isomorphic to
the subalgebra in $\tilde{V}_q^-\otimes\tilde{V}_q^+$ generated
by $K_\gamma^{-1}\otimes K_\gamma$,$\gamma\in P(\pi)$ and 
$U^-_q\otimes U^+_q$. The last algebra is a quantum solvable algebra
with Q1-Q4. The Theorem 3.4 and its Corollaries may be applied for
$\OG$ (see Cor.3.9)\\
{\bf Example 5}. One can consider some quantized
representations of $sl_q(2)$ and their corresponding motion
algebras [P2]. The algebra
$sl_q$ and some motion algebras are a quantum solvable algebras.
But they are not an iterated $q$-skew extension.
Restricting these motion algebras [P2] on Borel(or nilpotent) part, 
one can get new
examples quantum solvable algebras satisfying Q1-Q4. 

\section{Prime factors are pure quantum}

The goal of this section is to prove that prime factors of quantum solvable
algebras, obeying conditions Q1-Q4, are pure quantum (see Proposition 2.9
and Remark after it).\\
{\bf Proposition 2.1}. 
Let $R$ be a Noetherian domain and a $C$-algebra generated by
$x_1,\ldots,x_n,x_{n+1}^{\pm 1},\ldots,x_{n+m}^{\pm 1}$.
Suppose that $R$ admits the diagonal group $H$
of automorphisms (i.e $h(x_i)=\lambda_i(h)x_i$, $\lambda_i(h)\in C^*$).
Suppose that a group $\Gamma$ generated by $\lambda_i(h)$, $h\in H$ is
torsionfree
(in particular, all conditions is true if $R$ is
a quantum solvable algebra obeying Conditions Q1-Q3).\\
1) Let $\II$ be a prime ideal in $R$.
Consider maximal $H$-stable ideal $\II_H$ in $\II$.
That is 
$$\II_H=\bigcap_{h\in H}h(\II).$$
We assert that $\II_H$ is also prime.\\
2) Let $\PP$ be an $H$-stable semiprime ideal in $R$ and 
$\QQ$ is minimal prime ideal over $\PP$,
then $\QQ$ is $H$-stable.\\   
{\bf Proof.}
1)
The ideal $\II_H$ is always $H$-prime ([GL1], 2.1-2),
in particular $\II_H$ is semiprime. 
Denote by 
$\Lambda=\{\QQ_1,\ldots,\QQ_s\}$ the set of minimal prime ideals
containing $\II_H$. We have $\II_H=\QQ_1\bigcap\ldots\bigcap\QQ_s$.
The group $H$ acts on $\Lambda$ by permutations and there is only one 
$H$-orbit in $\Lambda$. Denote $H_1$ the stabilizer of $\QQ_1$ in $H$.
The subgroup $H_1$ is a subgroup of finite index in $H$.
Since $H$ is a commutative finitely generated group,
there exists a positive integer $N$ such that 
the subgroup $H^N=<\tau_1^N,\ldots,\tau_{n+m}^N>$ is contained in $H_1$. 

We can decompose an arbitrary element $a=a_1+\ldots+a_m$ into a sum 
of $H$-weight components $a_i$.
That is $h(a_i)=\chi_i(h)a_i$, 
$\chi_i(h)=\prod \lambda_i^{s_{i}(h)}$ with $s_{i}(h)\in \N$ for 
$i\in\overline{1,n}$ and 
$s_{i}(h)\in \Z$ for 
$i\in\overline{n+1,n+m}$.
Here the characters $\chi_i$, $\chi_j$ are different, if $i\ne j$.
For every $h$ the weight $\chi_i(h)$ is an element
of the group $\Gamma$.

Let us compare $H^N$-weights of $a_1,\ldots,a_m$.
The $H^N$-weight of $a_i$ equals to
$\chi_i(h^N)=(\chi_i(h))^N$. If $\chi_i(h^N)=\chi_j(h^N)$ for all $h$,
then $(\chi_i(h))^N=(\chi_j(h))^N$. The group $\Gamma$ is torsionfree.
It implies that $\chi_i(h)=\chi_j(h)$ for all $h$. Therefore $i=j$.
We prove that the $a_1,\ldots,a_m$ have different
$H^N$-weights.

Suppose that $a\in\QQ_1$. The ideal $\QQ_1$ is $H^N$-stable.
Therefore all $H^N$-components of $a$ (i.e $a_i$) also lie in $\QQ_1$.
But $\{a_i\}$ is indeed $H$-components of $a$.
We prove that, if $a\in\QQ_1$, then all its $H$-weight components lie in 
$\QQ_1$. This proves that $\QQ_1$ is $H$-stable.
There is only one $H$-orbite on $\Lambda$. Then $\Lambda=\{\QQ_1\}$
and $\II_H=\QQ_1$. It follows that $\II_H$ is prime ideal.\\ 
2) The $Q$ contains
$\PP$ and $\PP$ is $H$-stable. We have
$\QQ\supset\QQ_H\supset\PP$. We have proved that $\QQ_H$ is prime.
Since $\QQ$ is minimal over $\PP$ we have $\QQ=\QQ_H$.$\Box$\\ 
{\bf Lemma 2.2}. Let $B'=B[x;\tau,\delta]$ be an skew extension.
Let $S$ be a $H$-stable Ore subset in $B$.
Then $S$ is an Ore subset in $B'$.\\
{\bf Proof}. The subset $S$ is multiplicatively closed subset
 of regular elements in $B$ and $B'$.
It is sufficient to prove that for every element $s\in S$ and positive integer
$n$ there
exists $t,t'\in S$ and $b,b'\in B'$ such that $tx^n=bs$ and 
$x^nt'=sb'$. We shall prove the existence of $t,b$. The existence
of $t',b'$ is proved similaly. We shall prove this by induction on $n$.

First put $n=1$. Recall that for $xa=\tau(a)x+\delta(a)$.
We get $\tau(s)x=xs-\delta(s)$. For two elements
$\delta(s)\in B$, $s\in S$ there exist
the elements
$s_1\in S$ , $b_1\in B'$ such that $s_1\delta(s)=b_1s$.
We get $s_1\tau(s)x=s_1xs-s_1\delta(s)=
s_1xs-b_1s=(s_1x-b_1)s$. Denoting $t=s_1\tau(s)\in S$ and
$b=s_1x-b_1$, we complete case $n=1$.

Suppose that the statement is true for $n-1$. We are going to prove for $n$.
We have proved the existence of $t\in S$, $b\in B'$ such that
$tx=bs$. By assumption of induction, there exist
$t^*\in S$ and $b^*\in B'$ such that
$t^*x^{n-1}=b^*t$. We get $t^*x^n=(t^*x^{n-1})x=b^*tx=b^*bs$.
This completes the proof.$\Box$\\ 
We recall that $\delta$ is an inner $\tau$-derivation of a ring $R$, if
there exists $r\in R$ with the property
 $\delta(a)=\tau(a)r-ra$ for all $a\in R$.\\
{\bf Lemma 2.3}. Let $A$ be a division algebra over a field $F$,
$\mbox{char}(F)=0$. Let $A'=A[x;\tau,\delta]$ be a $q$-skew extension
(i.e $\delta\tau=q\tau\delta$ with $q\in F^*$).
Suppose that $q$ is not a root of unity and $\delta$ is 
not an inner $\tau$-derivation of $A$.
Then the ring $A'$ is simple.\\
{\bf Proof}. Well known (GL1,2.5) that for every $a\in A$ there holds
$$x^na=\sum_{i=0}^n{n\choose i}_q\tau^{n-i}\delta^i(a)x^{n-i}=
\tau^n(a)x^n+{n\choose 1}_q\tau^{n-1}\delta(a)x^{n-1}+\ldots+\delta^n(a).$$
Here ${n\choose i}_q=\frac{(n)!}{(i)!(n-i)!}$ where 
$(n)=\frac{q^n-1}{q-1}$ in the case $q\ne 1$ and $(n)=n$ in the case $q=1$.
It follows that
$$x^n\tau^{-n}(a)=ax^n+{n\choose 1}_q\tau^{n-1}\delta(\tau^{-n}(a))x^{n-1}+
\ldots +\{\mbox{terms of lower degree}\} $$
$$= ax^n+{n\choose 1}_qq^{-(n-1)}\delta(\tau^{-1}(a))x^{n-1}+\ldots
\{\mbox{terms of lower degree}\}.$$
Denote $\gamma={n\choose 1}_qq^{-(n-1)}=(n)q^{-(n-1)}\ne 0$.
We get 
$$
x^n\tau^{-n}(a)= ax^n+\gamma\delta(\tau^{-1}(a))x^{n-1}+\ldots
\{\mbox{terms of lower degree}\}.$$

Let $\II$ be an non zero ideal in $A'$. Then $\II=A'f(x)$ where
$f(x)=x^n+c_1x^{n-1}+\ldots+c_n$, $c_i\in A$ is the
 polynomial of the lowest degree in $\II$.
For every $a\in A$ consider the element 
$b=f(x)\tau^{-n}(a)-af(x)\in\II$.
By direct calculations, we obtain
$$b=f(x)\tau^{-n}(a)-af(x)=(x^n+c_1x^{n-1}+\ldots+c_n)\tau^{-n}(a)
-a(x^n+c_1x^{n-1}+\ldots+c_n)$$
$$
=
[\gamma\delta(\tau^{-1}(a))+c_1\tau^{-1}(a)-ac_1]x^{n-1}+
\{\mbox{terms of lower degree}\}.$$
Since the degree of $f(x)$ is lowest among polynomials in $\II$,
$\gamma\delta(\tau^{-1}(a))+c_1\tau^{-1}(a)-ac_1=0$.
Substituting $a$ by $\tau(a)$ and
denoting $r=\frac{c_1}{\gamma}$, we get $\delta(a)=\tau(a)r-ra$.
We get a contrudiction.
The algebra $A'$ is simple.$\Box$\\   

 Let $B$ be a $C$-algebra (i.e. $C$ is contained in the center of $B$).
We consider speciallization  $B_\eps=B/B\EE$.
For every element $a\in B$ we denote by $a_\eps$ its image under the map
$B\mapsto B_\eps$.\\
{\bf Definition 2.4}. Let $B$ be a Noetherian domain and a $C$-algebra.
We say that the $\EE$-Property holds for $B$, if\\ 
1) there exist a nonempty Zariski-open set $\Lambda$ of specializations
such that $B_\eps$ is semiprime for all $\eps\in\Lambda$,\\
2) for any nonzero element
$b\in B$ there exists a nonempty  
Zariski-open set $\Lambda_b$ of specializations
such that $b$ is a regular element in $B_\eps$.

Further we say that a property is true for
almost all specializations $\eps$ if this property is true for some
nonempty Zariski-open subset of prime ideals $\EE$ in $\Spec(C)$.

For a Noetherian domain $R$ we consider $S=B^*$, $A=Fract(B)=BS^{-1}$.
If $B$ obeys the $\EE$-Property, we denote by $S_\eps$,$\eps\in\Lambda$ 
the set of regular elements in $B_\eps$ and by $A_\eps$ the Goldie
quotient ring of $B_\eps S_\eps^{-1}$.
For any element $y=a^{-1}b\in A$, 
$a\in S$,$b\in B$ there exists a nonempty Zariski-open set $\Lambda$ of specializations of $C$ such that $a_\eps$ is 
regular in $B_\eps$ and $B_\eps$ is semiprime. 
Consider the element 
$y_\eps=a_\eps^{-1}b_\eps$ in the ring 
$A_\eps$.
We shall refer the element $y_\eps$ as the specialization of $y$
and prove latter that $y_\eps$ does not depends on presentation $y$.

Let $B'=B[x;\tau,\delta]$ be an Ore extension of 
$B$ with $\EE$-Property. 
By Lemma 2.2, $S$ is an Ore subset in $B'$.
We denote by $A'$ the ring $A[x;\tau,\delta]=B'S^{-1}$. 
By Lamma 2.2, $S_\eps$ is an Ore subset in $B'_\eps =B_\eps[x,\tau,\delta]$.
We consider  and $A'_\eps=B_\eps S_\eps^{-1}=A_\eps[x;\tau,\delta]$
Similar to above we can define the specialization of any
element $y=a^{-1}b$,$a\in S$, $b\in B'$.\\
{\bf Lemma 2.5}. For any element $y\in A'$ 
the element $y_\eps$ does not depends
on the presentation $y=a^{-1}b$,$a\in S$, $b\in B'$ for almost all $\eps$.\\
{\bf Proof}. Let $y=a^{-1}b$,$a\in S$, $b\in B'$.\\
1) There exist the elements $d\in S$ and $c\in B'$ such that 
$y=a^{-1}b=cd^{-1}$. Then $ac=bd$. Hence $a_\eps c_\eps=b_\eps d_\eps$
and $d_\eps$,$a_\eps$ are regular  for
almost all $\eps$. Whence $a_\eps^{-1}b_\eps=c_\eps d_\eps^{-1}$.\\
2) If there exists two presentations $y=a^{-1}b=a_1^{-1}b_1$, then
$a_\eps^{-1}b_\eps=c_\eps d_\eps^{-1}$ and 
$a_{1\eps}^{-1}b_{1\eps}=c_\eps d_\eps^{-1}$. This proves the statement.
$\Box$\\
{\bf Lemma 2.6}. 1) For all $y,z\in A'$ there hold
$(y+z)_\eps=y_\eps+z_\eps$, $(yz)_\eps=y_\eps z_\eps$.\\
2) Let $f$ be an element of the free algebra 
$A\{t_1,\ldots,t_n\}$. For  $y_i\in A'$,
$i\in \overline{1,n}$  
we denote by $f\{y_1,\ldots,y_n\}\in A$ the result
of substitution $t_i\mapsto y_i$.
Let $f_\eps$ denote the specializaton of $f$ in 
$A_\eps\{t_1,\ldots,t_n\}$
and $y_{i\eps}$ the
specialization of $y_i$. Then  $f\{y_1,\ldots,y_n\}=0$ implies
$f_\eps\{y_{1\eps},\ldots,y_{n\eps}\}=0$ for almost all $\eps$.\\
{\bf Proof}. Point 2) can be derived from Point 1).
We shall prove the first property
$(y+z)_\eps=y_\eps+z_\eps$ for $y=a^{-1}b$, $z=c^{-1}d$.
There exist $c_1\in S$, $a_1\in B'$ such that $ac^{-1}=c_1^{-1}a_1$
(i.e. $c_1a=a_1c$). 
We have 
$$y+z=a^{-1}b+c^{-1}d=a^{-1}(b+ac^{-1}d)=a^{-1}(b+c_1^{-1}a_1d)=
(c_1a)^{-1}(c_1b+a_1d).$$
For almost all
$\eps$ we get 
$$(y+z)_\eps=((c_1a)_\eps)^{-1}(c_1b+a_1d)_\eps=
= (c_{1\eps}a_\eps)^{-1}(c_{1\eps}b_\eps+a_{1\eps}d_\eps),$$
$ c_{1\eps}a_\eps=a_{1\eps}c_\eps$ with regular elements
$a_\eps$, $c_\eps$, $c_{1\eps}$.
This implies
$$(y+z)_\eps=a_\eps^{-1}b_\eps+c_\eps^{-1}d_\eps= y_\eps+z_\eps$$
for almost all $\eps$.$\Box$\\

Let $\PP$ be an ideal in $R$.
Suppose that $\PP\bigcap C= 0$. For any specialization 
$\eps$ we consider the ideal
$$\PP_\eps=\frac{\PP+R\EE}{R\EE}$$ in $R_\eps.$
One can consider the map
$\pi_\eps:R \mapsto R/\PP\mapsto (R/\PP)_\eps$ where
$$
(R/\PP)_\eps=
=\frac{R/\PP}{(R/\PP)\EE}=\frac{R}{\PP+R\EE}=
\frac{R/{R\EE}}{(\PP+R\EE)/{R\EE}}= R_\eps/\PP_\eps.\eqno(2.1)$$
{\bf Proposition 2.7}. Let a quantum solvable algebra $R$ obeys
conditions $Q1-Q4$. Let $\PP$ be an $H$-stable prime ideal in
$R$. Suppose that $\PP\bigcap C=0$. We assert that the quotient ring
$B=R/\PP$ obeys $\EE$-Property.\\
{\bf Proof.}
We are going to prove the following two statements:\\
1) For almost all $\eps$ the ideal $\PP_\eps$ is semiprime
(in particular $\PP_\eps\ne R_\eps$);\\
2) Let $b\notin\PP$. For almost all $\eps$
 the image $\pi_\eps(b)$ is a regular element 
in $R_\eps/\PP_\eps$.
 We shall prove the Proposition by induction on filtration
$$R=R_1\supset\cdots\supset 
C_q[k_1^{\pm 1},\ldots,k_m^{\pm 1}]\supset\ldots\supset 
C_q[k_m^\pm]\supset C \supset 0.$$
First we notice that the assertion is evidently true for $R=C$. 
We assume that the statement is true for quantum solvable $R_1$ generated
by $x_1,\ldots,x_n,k_1^{\pm 1},\ldots,k_s^{\pm 1}$.
Consider a $q$-skew extension
$R'=R[x;\tau,\delta]$ satisfying the requires of the Proposition.
We are going to prove that the statement is true for $R'$
with prime ideal $\PP'$ obeying the assumptions of the Proposition.

We denote $\PP=\PP'\bigcap R$. Since the ideal $\PP'$ is prime, it is
completely prime (Remark after Condition Q2). 
It follows that $\PP$ is also completely prime.
By assuption, $\PP'$ is stable under the action of group
$H'=<\tau,\tau_1,\ldots,\tau_n>$.
Then $\PP$ is stable under the action of group
$H=<\tau_1,\ldots,\tau_n>$.
The ideal $\PP$ satisfies the conditions of the Proposition.
By assuption of induction, the Proposition is true for $\PP$.
The ideal $\PP$ is $\tau$-stable. Then $\PP$ is $\delta$-stable.
Hence $R'\PP$ is the two-sided ideal in $R'$.

We shall use the notations
$B=R/\PP$, $S=B^*= B-\{0\}$, $A= BS^{-1}=\Fract(B)$, 
$B'= R'/(R'\PP) = B[x;\tau,\delta]$.
As in Lemma 2.5 we consider $A'= A[x;\tau,\delta]=B'S^{-1}$.
Denote by $\Po$ an ideal $\PP'/(R'\PP)$ in $B'$. It is obvious that
$R'/\PP'=B'/\overline{\PP}$.

The statement is true in the case $\PP'=R'\PP$, by the assumption of induction.
Let us consider $\PP'\ne R'\PP$. 
 Then $A'$ is not a simple ring. According to Lemma 2.3 
 $\delta$ is an inner 
$\tau$-derivation of $A$ (i.e. there exists $r\in A$ such that 
$\delta(a)=\tau(a)r-ra$ for all $a\in A$). 
Denote $\xt=x+r$. Then $A'=A[\xt,\tau].$ The proof divides in two cases
$\xt\in \Po S^{-1}$ and $\xt\notin \Po S^{-1}$. We begin 
with the more complicated second case.\\
i) Suppose that $\xt\notin \Po S^{-1}$.\\ 
One knows the classification of prime ideals in $A'$ 
([J], Theorems 0,1; [GL1],2.3 and Theorem 11.1).
If $\tau^N$ is not an inner automorphism in $A$ for all $N$, then
$\PP'=R'\PP$ and the statement is true.
Suppose that for some $N$ the automorphism 
$\tau^N$ is inner , 
i.e. there exists $g\in A$ such that
$\tau^N(a)=g^{-1}ag$ for all $a\in A$.
Suppose that $N$ is minimal with this property.
We consider the element $z=g\xt^N$ commuting with $A$.  
The ideal $\Po S^{-1}$ have a form $A'f(z)$.
Here $f(t)=t^n+a_1t^{n-1}+\ldots+a_n$ is an irreducible polynomial over 
the field af $\tau$-invariants $k=(\Center(A))^{\tau}$.
Since ${\mbox{char}}(k)=0$, then $\Disc(f)\ne 0$.

The ideal $\Po$ is generated as a left ideal 
 by the finite number of elements
$\phi_1,\cdots,\phi_s\in B' = B[x,\tau,\delta]$
(i.e. $\Po=B'\phi_1+\ldots+B'\phi_s$). 
There exist the polynomials 
$$h_i(t)=\beta_0^{(i)}t^{n_i}+\ldots+\beta_{n_i}^{(i)}\quad 
i\in \overline{1,s}$$
in $A[t]$ such that
$$\phi_i=h_i(\xt)f(z).\eqno(2.2)$$
There also exist the polynomials
$$G_i=\alpha_0^{(i)}t^{m_i}+\ldots+\alpha_{m_i}^{(i)},\quad i\in
\overline{1,s}$$
in $A[t]$ such that 
$$G_1(\xt)\phi_1+ \ldots + G_s(\xt)\phi_s= f(z).\eqno(2.3)$$
All elements $$r, g, a_1,\ldots,a_n, \Disc(f),
\{\alpha_l^{(i)},\beta_j^{(i)}\}_{i\in\overline{1,s};j\in\overline{0,n_i}
,l\in\overline{0,m_i}}\eqno(2.4)$$
can be presented in the form $a^{-1}b$ with $a\in S$,$b\in B$ :
$$ r=(r_1)^{-1}r_2,\quad g=(g_1)^{-1}g_2\ne 0,
\quad\Disc(f)=(d_1)^{-1}d_2\ne 0, \eqno(2.5)$$ 
$$a_k=a_{k1}^{-1}a_{k2}, k\in\overline{1,n}\eqno(2.6)$$
$$\alpha_l^{(i)}=(\alpha_{l1}^{(i)})^{-1}\alpha_{l2}^{(i)}).\eqno(2.7)$$
$$\beta_j^{(i)}=(\beta_{j1}^{(i)})^{-1}\beta_{j2}^{(i)}.\eqno(2.8)$$

Denote by $\pht_1,\ldots,\pht_s$ the representatives of
$\phi_1,\ldots,\phi_s\in B'$ in $R'=R[x,\tau,\delta]$.
For any specialization $\eps$ we have
$$ 
\PP'_\eps=\frac{\PP'+R'\EE}{R'\EE}=
\frac{R'\pht_1+\cdots+R'\pht_s+R'\PP+R'\EE}{R'\EE}.\eqno(2.9)
$$
The left ideal $R'\PP+R'\EE$ consists
of polynomials in $R'=R[x,\tau,\delta]$ with coefficients in
$\PP+R\EE$. Then $(R'\PP+R'\EE)/R'\EE$ are polynomials in
$R'_\eps=R_\eps[x,\tau,\delta]$ with coefficients in 
$$\frac{\PP+R\EE}{R\EE}=\PP_{1\eps}.$$
We obtain
$$\frac{R'\PP+R'\EE}{R'\EE}= R'_\eps\PP_{\eps}$$
Consider 
$$B'_\eps=\frac{R'_\eps}{R'_\eps\PP_{\eps}}
\quad B_{\eps}=\frac{R_{\eps}}{\PP_{\eps}}$$
We get
$$B'_\eps=\frac{R'/R'\EE}{(R'\PP+R'\EE)/R'\EE}
=\frac{R'}{R'\PP+R'\EE}=\frac{R}{\PP+R\EE}[x;\tau,\delta]=
B_{\eps}[x;\tau.\delta] $$
Denote by $\Po_\eps$ the ideal 
$$\frac{\PP_\eps}{R_\eps\PP_{1\eps}}$$
in $B'_\eps$.
Denote by $\phi_{i\eps}$ the images
of $\pht_i$ in 
$B'_\eps=B_{\eps}[x;\tau,\delta]$.
Using (2.9), we obtain
$$\Po_\eps=
\frac{\PP'_\eps}{R'_\eps\PP_{\eps}}=
\frac{\frac{R'\pht_1+\cdots+R'\pht_s+R'\PP+R'\EE}{R'\EE}}
{\frac{R'\PP+R'\EE}{R'\EE}}$$
$$ = \frac{R'\pht_1+\cdots+R'\pht_s+R'\PP+R'\EE}{R'\PP+R'\EE}
=B'_\eps\phi_{1\eps}+\cdots+B'_\eps\phi_{s\eps}
$$
We denote by $r_{1\eps}$, $g_{1\eps}$, $g_{2\eps}$,
$d_{1\eps}$, $d_{2\eps}$, $a_{k1\eps}$, 
$\alpha_{l1\eps}^{(i)}$, $(\beta_{j1\eps}^{(i)})$
the specializations of the elements
$r_1$, $g_1$, $g_2$,
$d_1$, $d_2$, $a_{k1}$, 
$\alpha_{l1}^{(i)}$, $\beta_{j1}^{(i)}$ in $B_\eps$.
By assumption of induction, there exists a nonempty 
Zariski-open set $\Lambda$ of 
specialization such that the elements
$r_{1\eps}$, $g_{1\eps}$, $g_{2\eps}$,
$d_{1\eps}$, $d_{2\eps}$, $a_{k1\eps}$, 
$\alpha_{l1\eps}^{(i)}$, $(\beta_{j1\eps}^{(i)})$
 are regular in 
$B_{\eps}$ and specializations of
(2.2), (2.3) hold in $B'_\eps$ (by Lemma 2.6).

The Goldie quotient ring $A_{\eps}=B_{\eps}S_\eps^{-1}$ 
is the sum of simple algebras
$$A_{\eps}= \oplus_{i=1}^k Mat(D_i)$$
where $D_i$ are division rings.
The set of regular elements $S_\eps$ is an denominator set in 
$B'_\eps=B_{1\eps}[x,\tau,\delta]$ (see Lemma 2.2).
As above we denote ${A'_\eps}= B'_\eps S_\eps^{-1}=
A_{\eps}[x,\tau,\delta]$.
Using the decompositions (2.5-8) one can present the images
of elements (2.4) in $B_{\eps}$ for $\eps\in\Lambda$.
We denote these images by $r_\eps,g_\eps, d_\eps=\Disc(f)_\eps,\ldots$.
Recall that the elements $g_\eps$ and $d_\eps$ are invertible
in $A$. 

We denote $\xt_\eps=x+r_\eps\in A'_\eps$ and 
$z_\eps=g_\eps\xt_\eps^N$. 
Then $A'_\eps=A_{\eps}[\xt_\eps,\tau]$.
We denote by $G_{i\eps}(t)$, $f_\eps(t)$ the polynomials
with the correspoding coefficients in $A_{\eps}$.
We get 
$$\overline{\PP}_\eps S_\eps^{-1}=A'_\eps\phi_{1\eps}+\cdots+
A'_\eps\phi_{s\eps}=
A'_\eps f_\eps(z_\eps)$$
In particullar, $\overline{\PP}_\eps S_\eps^{-1}\ne A'_\eps$.

The coefficients $a_{k\eps}$ of $f_\eps$ 
lie in $k_\eps= (\Center A_{\eps})^\tau$.
The automorphism $\tau$ acts on $A_{\eps}$.
Here $k_\eps$ is not
a field in general but a sum $k_\eps=K_1\oplus\cdots\oplus K_M$
of the fields $K_i$.
The algebra $A_{\eps}$ decomposes into the sum 
$A_{\eps}=A_1\oplus\cdots\oplus A_M$ of $\tau$-prime algebras.
Consider the system $e_i, i\in\overline{1,M}$ 
of minimal $\tau$-invariant central idempotents in $A_\eps$.
Since $A'_\eps =A_\eps[\xt_\eps;\tau]$, the 
elements $e_i$ lie the center of $A'_\eps$ 
We have  $A_i=A_{\eps}e_i$ and  
we get the decomposition 
$A'_\eps=A'_\eps e_1\oplus\ldots\oplus A'_\eps e_M.$
We denote further $A'_{i\eps}=A'_\eps e_i$.
The polynomial $f_\eps(t)e_i$ decomposes over $K_i$ into product of
$K_i$-irreducible polynomials.
Since $d_\eps=\Disc(f_\eps(t))$ is an invertible element in $A_{\eps}$,
the polynomial $f_\eps(t)$ have different and nonzero 
(since $\xt\notin\Po S^{-1}$) roots in  
algebraic closure of the field  $K_i$.
It follows that multiplicites of irreducible polynomials are equal to one.
We have  
$f_\eps(t)e_i=\prod p_{ij}(t)$.

Accoding to (2.1), 
$$(R'/\PP')_\eps=R'_\eps/\PP'_\eps=B'_\eps/\overline{\PP}_\eps.$$
We get decomposition
$$
\left(\frac{R'_\eps}{\PP'_\eps}\right)_{S_\eps}= 
\frac{B'_\eps S_\eps^{-1}}{\overline{\PP}_\eps S_\eps^{-1}}
=\frac{A'_\eps}{A'_\eps f_\eps(z_\eps)}=
\oplus_{i=1}^M 
\frac{A'_{i\eps}}{A'_{i\eps}f_\eps(z_\eps)}=
\oplus_{ij}\frac{A'_{i\eps}}{A'_{i\eps}
p_{ij}(z_\eps)}.
$$
Since $\tau^N(a)=g_\eps^{-1}ag_\eps$ in $A_{\eps}$, the automorphism 
$\tau^N$ in inner in $A_{\eps}$.
Let $N_1$ be the smallest positive integer
such that $\tau^{N_1}$ is inner in $A_{\eps}$. 
 Then $\tau^{N_1}(a)=h_\eps^{-1}ah_\eps$
for some invertible elements $h_\eps\in A_{\eps}$.
Denote $y_\eps=h_\eps\xt_\eps^{N_1}$.
The integer $N_1$
divides $N$. Denote $L=\frac{N}{N_1}$.
There exists central invertible element $\gamma\in A_{\eps}$
such that 
$z_\eps=\gamma y_\eps^L$.
We get
$$ \left(\frac{R'_\eps}{\PP'_\eps}\right)_{S_\eps}= 
\oplus_{ij}\frac{A'_{i\eps}}{A'_{i\eps} 
p_{ij}(\gamma y_\eps^L)}.
$$
As we notice above the polynomial $f_\eps(t)$ have 
different and nonzero roots in algebraic closure of every $K_i$,
the roots of polynomial $f_\eps(\gamma t^L)$ are also different.
We get decompositions into products of irreducible polynomials
with multiplicities equal to one
$p_{ij}(t)=\prod p_{ijk}(t)$.
Finally,
$$ \left(\frac{R'_\eps}{\PP'_\eps}\right)_{S_\eps}= 
\oplus_{ijk}
\frac{A'_{i\eps}}{A'_{i\eps}p_{ijk}(y_\eps)}.
$$
This proves of 
that $R'_\eps/\PP'_\eps$ is semiprime (i.e. statement 1) of Proposition)
in the case i).

Let us prove statement 2) in the case i) .
Consider the case $\PP'\ne R'\PP$, $\xt\notin\PP'$.
As above the ideal $(\PP'/R'\PP)S^{-1}$ equals to $A'f(z)$. 

Consider the element $b\in R'- \PP'$.
Let $\bar{b}$ is the image of $b$ in $B'=R'/\PP'$.
 We treat $\bar{b}$ as polynomial $\bar{b}=\bar{b}(\xt)$ in $A'=A_1[\xt,\tau]$.
Let $d(\xt)$ be left GCD ([B],VIII,12,Ex.7)  
of $\bar{b}(\xt)$ and $f(z)=f(g\xt^N)$.
Since $\PP'$ is completely prime, 
the polynomial  $f(z)$ can not 
decompose into product of polynomials of lower 
degree. It follows that $d(\xt)=d\in A$.

There exist the polynomials $M_1(\xt)$ and $M_2(\xt)$ such that
$$M_1(\xt)\bar{b}(\xt)+M_2(\xt)f(g\xt^N)=d.$$
Similar to the proof of statement 1), for almost all $\eps$
the ideal $\PP'_\eps$ is semiprime and 
$d$ can be represented as an invertible element in $A_{\eps}$.
By Lemma 2.6, we 
get
$$M_{1\eps}(\xt_\eps)\bar{b}_\eps+M_{2\eps}(\xt_\eps)f(z_\eps)=d_\eps.$$
It implies
$$M_{1\eps}(\xt_\eps)\bar{b}_\eps = d_\eps\mod(\Po S^{-1}).$$
We see that $\bar{b_\eps}$ is not left zero divisor.
Since $\PP'_\eps$ is semiprime then $\bar{b_\eps}$ is not a right
zero divisor.We complete the proof in the case i).\\  
ii) Suppose that $\xt\in \Po S^{-1}$. 
The proof in 
this case is similar to case i), putting $f(t)=t$ and $z=\xt$. $\Box$\\
{\bf Definition 2.8}.([P1]) Let $B$ be a Noetherian domain.
We say that $B$ is pure quantum if 
there is no elements $x,y$ in $\Fract(B)$ related as follows 
$xy-yx=1$.\\
Recall that Weyl algebra $\Bbb{A}_1(F)$ is an algebra over a field $F$ 
and generated by elements
$X,Y$ with the single relation  $XY-YX=1$. Well known that 
$\Bbb{A}_1(F)$  is simple 
over a field of zero characteristic. 
Therefore an $F$-algebra $B$ is pure quantum if there exist
an 
embedding of the Weyl algebra $\Bbb{A}_1(F)$ in $\Fract(B)$
([P1]).\\
{\bf Proposion 2.9}. Let $\PP$ be an $H$-stable prime ideal in 
quantum solvable algebra $R$,
satisfying Q1-Q4 and $\PP\bigcap C=0$. Then the ring $B=R/\PP$ 
is pure quantum.\\
{\bf Proof}.
Let $x,y$ be an elements with $xy-yx=1$ in $\Fract(B)$. 
Accoding to Lemma 2.6 and Proposition 2.7, 
there exists a nonempty Zariski-open set $\Lambda$ in 
$\Spec(C)$ such that 
$x_\eps y_\eps - y_\eps x_\eps=1$ for $\EE=\Ker(\eps)\in \Lambda$ .  
Here $x_\eps,y_\eps\in A_\eps=\Goldie(B_\eps)$.
On other hand  there exists a dense set $\Omega\in\Spec(C)$ 
(Condition Q4) for which $A_\eps$ is finite over its 
center $F_\eps$.
There exists a common point $\eps_0\in\Lambda\bigcap\Omega$.
 Consider the homomorphism $X,Y\mapsto x_{\eps_0},y_{\eps_0}$ 
of the algebras
$\Bbb{A}_1(F_{\eps_0})\mapsto A_{\eps_0}$.
The first algebra is infinite
over $F_{\eps_0}$ and the second is finite over the same field.
There exists non zero kernel. But the first algebra is simple. 
This leads to contradiction. The ring $B$ is pure quantum.
$\Box$\\
{\bf Remark}. We shall prove in the next section that the claim 
of $H$-stablity of $\PP$ can be omited here (Corollary 3.7).   

\section{Stratification of spectrum}

In this section we study 
stratification of prime ideals (Theorem 3.4, Theorem 3.10, Example 3.12). 
This gives
stratification of prime ideals over $F=\Fract(C)$ (Corollary 3.9)
and proves conjecture on skew field of fractions (Corollary 3.8).\\
{\bf Definition 3.1}. Let $P=(p_{ij})_{ij=1}^m$ be a matrix with 
entries $p_{ij}\in C^*$ (see Section 1) and  
$p_{ij}p_{ji}=p_{ii}=1$. An algebra 
$C_P[Y_1^{\pm 1},\ldots,Y_m^{\pm 1}]$ of twisted Laurent 
polynomials is generated by
$Y_1^{\pm 1},\ldots,Y_m^{\pm 1}$ related as follows
$Y_iY_j=p_{ij}Y_jY_i$.\\
{\bf Notation 3.2.} We denote by $N$ the multiplicatively closed set in $C$ 
generated by $\gamma-\gamma'$ where $\gamma,\gamma'\in\Gamma$ and
$\gamma\ne\gamma'$.\\
{\bf Lemma 3.3} Let $B$ be a $C$-algebra and a Noetherian semiprime ring
with the set $\{\QQ_1,\ldots,\QQ_s\}$ of minimal prime ideals.
Suppose that all $\QQ_j$ are completely prime and 
$B/_{Q_j}$ is pure quantum.

Suppose that $B$ is generated by the elements
$x_1, x_2,\ldots, x_n, k_1^{\pm},\ldots,
k_m^{\pm}$, $x_n\ne 0$ with the following relations. 
The elements $k_1,\ldots,k_m$ $q$-commute with all generators.
For $i<j$ the holds 
$x_ix_j=q_{ij}x_jx_i+r_{ij}$
where $r_{ij}$ is the element of subalgebra
$B_{i+1}$ generated by $x_{i+1},\ldots, x_n, k_1^{\pm},\ldots,
k_m^{\pm}$. 
Suppose that $x_n\notin \QQ_j$ for all $j$.
Suppose that $B$ admits $H$-action (see Condition Q3).
We assert the following.\\
1) The multiplicatively closed set generated by $x_n$ is a 
denominator set in $R$;
denote by $S$ the multiplicatively closed set generated by $x_n$ and $N$;\\
2) Let $B_i'=(B_i)_S$ be localization of $B_i$ over $S$.  
 There exist the new set of generators
$x'_1,\ldots,x'_{n-1}, x^{\pm 1}_n, k_1^{\pm},\ldots,
k_m^{\pm}$ of $B'$ such that  
all generators are $H$-weight vectors; $B'_i$ is generated by
$x'_i$ and $B'_{i+1}$; $x'_i\in B_i$ and  $x_n$  $q$-commutes with all new
generators and 
$x'_ix'_j=q_{ij}x'_jx'_i+r'_{ij}$
where $r'_{ij}\in B'_{i+1}$.\\
{\bf Remark 1}. Comparing with Definition 1.2 we observe that 
we don't claim that monomials form a $C$-basis. There may be some other
relations among generators besides (1.1).
The generators may even coinside and some of them may be equal to zero.\\
{\bf Proof}. The proof is similar to [P1].
The element $x=x_n$ is an elements of finite adjoint action if the following 
sense.
For every $y\in B$ there exists a polynomial
$f(t)= a_0t^k + a_1t^{k-1} +\cdots + a_k$, $a_0\ne 0$, $a_k\ne 0$ 
and $a_i\in C$ such that 
$ 
a_0x^ky + a_1x^{k-1}yx +\cdots + a_kyx^k = 0 $ 
and all roots of $f(t)$ belongs to $\Gamma$
([P1], Def.3.1,Lemma 4.3).

We shall prove that the element $x=x_n$ is regular in $B$.
Suppose that $ax=0$, $a\in B$.
Then $ax=0\mod(Q_j)$ for all $j$. Since $Q_j$ is completely prime
and $x\ne 0\mod(Q_j)$, we get $a=0\mod(Q_j)$ for all $j$. Whence $a=0$.
 
Moreover, the multiplicatively closed set generated by $x$ satisfies Ore
 conditions
([P1],Prop.3.3). Whence $S$ is a denominator set. This proves 1).

For all $j$ the algebra $B/_{Q_j}$ is pure quantum. 
This implies that the $\Ad_x$-action 
($\Ad_x(y)=xyx^{-1}$) in $B'$ is diagonalizable.
The proof is similar to ([P1];Prop.3.4).
If it is not true, there exist the nonzero 
elements $\overline{u},\overline{v}\in B'$
such that
$\Ad_x(\overline{u})=\alpha\overline{u}$,
$\Ad_x(\overline{v})=\overline{u}+\alpha\overline{v}$.
Since $\overline{u}\ne 0$ there exists $\QQ_j$ such that
$\overline{u}\ne 0\mod \QQ_j$. Denote $u=\overline{u}\mod(\QQ_j)=\ne 0$,
$v=\overline{v}\mod(\QQ_j)$. We have
$\Ad_x(u)=\alpha u$,
$\Ad_x(v)=u+\alpha v$. The ring $R/Q_j$ is a domain.
Consider $Y=\alpha u^{-1}vx^{-1}\in\Fract(R/\QQ_j)$ and check $xY-Yx=1$.
This contradics the claim that $R/Q_j$ is pure quantum.
We prove that
$\Ad_x$-action is diagonalizable.

Consider one of generators $x_i$. 
Suppose that $B'_i\ne B'_{i+1}$ (if there are equal, we put $x'_i=x'_{i+1}$).  
Let 
$\mu(t)=(t-\gamma_0)(t-\gamma_1)\cdots(t-\gamma_M)$, $\gamma_j\in\Gamma$ 
be the minimal polynimial of $\Ad_x$-action on $x_i$.
The generator $x_i$ can be decomposed into the sum of $\Ad_x$-eigenvectors
$x_i= x_i^{(0)}+ x_i^{(1)}+\cdots+x_i^{(M)}$ with 
$\Ad_x(x_i^{(m)})=\gamma_m x_i^{(m)}$.

We can calculate the eigenvectors directly
$$\mu_m(t)= (t-\gamma_0)\cdots\widehat{(t-\gamma_m)}\cdots(t-\gamma_M).$$
$$
x_i^{(m)}=\frac{1}    
{\prod_{j\ne m}(\gamma_m-\gamma_j)}\mu_m(Ad_x)x_i$$

Since $x_i$ and $x$ are $H$-weight vectors, the same is 
all $x_i^{(m)}$.
One of the eigenvalues, say $\gamma_0$, equals to $q_{ni}$.
All eigenvectors besides $x_i^{(0)}$ belong to $B'_{i+1}$,
but $x_i^{(0)}$ does not belongs to $B'_{i+1}$.
We put $x'_i=x_i^{(0)}$. We get new set of generators in $B'$ related 
as required ([P1], Main Theorem).

Finally, after right multiplication of $x'_i$ by some positive 
power of $x$ we get
the system of generators contained in $B$.
$\Box$\\
{\bf Remark 2}. Consider the finite set $N_1\in C$ of all denominators 
that occur in the above expression for $x^{(m)}_i$. Let $S_1$ be 
the denominator
set generated by $x$ and $N_1$. We notice that all $x'_i$ lie in
$(B')_{S_1}$. Further referring Lemma 3.3 we 
shall substitute further $S$ by finitely generated $S_1$.\\
{\bf Theorem 3.4}. Let $R$ be a quantum solvable algebra obeying Conditions
Q1-Q4.\\
1) There exists the finite set $\M=\{\PP_{i_1\ldots i_{k+1}}\}$
of $H$-stable semiprime ideals with zero intersection with
$C$ satisfying the following property.
There exists a denominator set $S_k$ 
generated by $k$ $q$-commuting elements in $R/\PP_{i_1,\ldots i_{k+1}}$,
such that 
$(R/\PP_{i_1,\ldots i_{k+1}})_{S_k}$ 
is isomorphic to a factor ring of the algebra
of twisted Laurent polynomials.
Further we refer an ideal from $\M$ as a standart ideal and
$S_k$ - as a standart denominator set.\\
2) Any prime ideal $\II$ of $R$ with zero intersection with $C$
contains a unique standart ideal
$\PP_{i_1,\ldots i_{k+1}}$ such that 
$S_k\bigcap(\II\mod(\PP_{i_1,\ldots i_{k+1}}))= \emptyset$.\\
{\bf Proof}. The proof will be devided into two steps.
In the first step we shall construct the $H$-stable  semiprime ideal
$\PP_{i_1,\ldots,i_{k+1}}$ with zero intesection with $C$.
Here $\{i_1,\ldots,i_{k+1}\}$ is the system of positive integes
obeying
$k+i_1+\cdots+i_{k+1}= n$. This will prove 1).
Indeed, the number of elements of $\M$ may be less then the number of 
nonnegative integer solutions of this equation.
The ideals of the form $ \PP_{i_1,\ldots i_{k+1}}$ may coinside and not
every solution $i_1,\ldots,i_{k+1}$  correspondes to standart ideal.
See the construction of $\PP_{i_1,\ldots i_{k+1}}$ and Remark 1 to Lemma 3.3.
The statement 2) will be proved in Step 1.\\
{\bf Step 1}. Consider the system of generators in reverse order
$x_n,x_{n-1},\ldots,x_1$.
Consider the set of the first $i_1$ elements in this system
$L_1=\{x_n,\ldots,x_{n-i_1+1}\}$.

 Denote by $\JJ=\JJ_{i_1}$ the minimal semiprime ideal which
contains the elements $L_1$.
For $i_1=0$ we put $L_1=\emptyset$ and $\JJ=0$.
We denote by $X_1=\{\QQ_1,\ldots,\QQ_s\}$. 
the set of minimal prime ideal contaning $\JJ$.
Since $\JJ$ is semiprime, $\JJ=\QQ_1\bigcap\cdots\bigcap\QQ_s$.
We decompose $X_1=X_1'\bigcup X_1''\bigcup X_1'''$
where
$$X_1'=\{\QQ\in X_1:\QQ\bigcap C= 0, x_{n-i_1}\notin\QQ, \},$$ 
$$X_1''=\{\QQ\in X_1:\QQ \bigcap C\ne 0, x_{n-i_1}\notin\QQ\},$$
$$X_1'''=\{\QQ\in X_1: x_{n-i_1}\in\QQ\}.$$ 
We denote
$$\PP_{i_1}=\bigcap_{Q\in X_1'}Q$$ 
The ideal $\PP_{i_1}$ contains $L_1$, does not contain 
$x_{n-i_1}$ and $\PP_{i_1}\bigcap C=0$.
We consider the semiprime ring $B_{i_1}=R/\PP_{i_1}$ generated by
$x_{n-i_1},\ldots,x_j,\ldots$ with $n-i_1\ge j$ and the homomorphism
$\pi_1:R\mapsto B_{i_1}$.
We save the notation $C$ for the image $\pi_1(C)$ in $B_{i_1}$.
The generators of $B_{i_1}$ satisfies the conditions of Lemma 3.3.

The ideal $\PP_{i_1}$ is $H$-stable and semiprime. 
All minimal prime ideals $Q\in X_1'$ are completely prime (Condition Q2).
They are $H$-stable (Lemma 2.1) and have zero intersection with $C$. 
Therefore, all domains
$R/\QQ$, $\QQ\in X_1'$ are pure quantum (Proposition 2.9).
 We have proved that $B_{i_1}$ satisfies all conditions of Lemma 3.3.
Denote by $S_1$ the denominator set in $B_{i_1}$ generated by $x_{n-i_1}$ and
finitely generated $N_1\subset C$ 
(see Notation 3.2, and Remark after Lemma 3.3).
Further we save the notations of generators of $B_{i_1}S_1^{-1}$
and denote $C_1=(C)_{N_1}$. 
By Lemma 3.3, the algebra $B_{i_1}S_1^{-1}$ has the new system of generators
$\{x_{n-i_1}^{\pm 1}, x'_{n-i_1-1},\ldots, x'_j\ldots\}$,
$x'_j\in B_{i_1}=\pi_1(R)$.

Consider the set of $i_2$ generators 
$L_2=\{x'_{n-i_1-1},\ldots,x'_{n-i_1-i_2}\}$.
Let $\JJ'_{i_1i_2}$ be the minimal semiprime ideal
in $B_{i_1}S_1^{-1}$ containing $L_2$.
We denote by $X_2$ the set of minimal prime ideals
in $B_{i_1}S_1^{-1}$ over $\JJ'_{i_1i_2}$.
As above we decompose
$X_2=X_2'\bigcup X_2''\bigcup X_2'''$
where
$$X_2'=\{\QQ'\in X_2:\QQ'\bigcap C_1= 0, x'_{n-i_1-i_2-1}\notin\QQ, \},$$ 
$$X_2''=\{\QQ'\in X_2:\QQ' \bigcap C_1\ne 0, x'_{n-i_1-i_2-1}\notin\QQ\},$$
$$X_2'''=\{\QQ'\in X_2: x'_{n-i_1-i_2-1}\in\QQ\}.$$
We denote
$$\PP'_{i_1i_2}=\bigcap_{Q'\in X_2'}Q'$$ 
As above   
$\PP'_{i_1i_2}$ an $H$-stable semiprime
ideal in $B_{i_1}S_1^{-1}$ such that it contains $L_2$, 
does not contain $x'_{n-i_1-i_2-1}$ and 
$\PP'_{i_1i_2}\bigcap C=0$.
We denote 
$$B_{i_1i_2}=\frac{B_{i_1}S_1^{-1}}{\PP'_{i_1i_2}}$$
$$\pi_2:B_{i_1}\mapsto B_{i_1}S^{-1}\mapsto B_{i_1i_2}$$
We save the notation $S_1$ to the image of $S_1\in B_{i_1}$ under $\pi_2$.

The algebra $B_{i_1i_2}$ is generated by the system of elements
$$\{x_{n-i_1}^{\pm 1},x'_{n-i_1-i_2-1},\ldots,x'_j,\ldots,\}.$$
We shall prove that $B_{i_1i_2}$ satisfies the conditions of Lemma 3.3. 

Recall that if an Noetherian domain $A$ admits a denominator set $S$, then
there exists one to one correspondence 
$\II\mapsto \JJ=\II\cap A$, $\JJ\mapsto\II=\JJ S^{-1}$ between
ideals $\II$ in $AS^{-1}$ and ideals $\JJ$ in $A$ with the property
that $as\in \JJ$ or $sa\in\JJ$ implies $a\in\JJ$ ([D],3.6.15).

According to this property,
$\PP'_{i_1i_2}\bigcap B_{i_1}$ is $H$-stable semiprime ideal in $B_{i_1}$
and $\{\QQ'\bigcap B_{i_1}: \QQ'\in X_2'\}$ is the set of all minimal
prime ideals over $\PP'_{i_1i_2}\bigcap B_{i_1}$ .
Denote
$$\PP_{i_1i_2}=\pi_1^{-1}(\PP'_{i_1i_2}\bigcap B_{i_1})=
\pi_1^{-1}\pi_2^{-1}(0).$$
We notice $\PP_{i_1i_2}\supset \PP_{i_1}$ and
the set 
$$\{\QQ=\pi_1^{-1}(\QQ'\bigcap B_{i_1}): \QQ'\in X_2'\}$$
is the set of all minimal prime ideals over
$\PP_{i_1i_2}$.
We have 
$$\frac{R}{\PP_{i_1i_2}}=\frac{R/P_{i_1}}{\PP_{i_1i_2}/\PP_{i_1}}=
\frac{B_{i_1}}{\PP'_{i_1i_2}\bigcap B_{i_1}},$$
$$\left(\frac{R}{\PP_{i_1i_2}}\right)_{S_{i_1}}=
\frac{B_{i_1}S_1^{-1}}{(\PP'_{i_1i_2}\bigcap B_{i_1})S_1^{-1}}
= \frac{B_{i_1}S_1^{-1}}{\PP'_{i_1i_2}}=B_{i_1i_2}.$$ 
The  minimal prime ideals in $B_{i_1i_2}$
are 
$\QQ'/\PP'_{i_1i_2}$ with $\QQ\in X'_2$.
We get
$$(R/\QQ)_{S_1}=\left(\frac{R}{\pi_1^{-1}(\QQ'\bigcap B_{i_1}}\right)_{S_1}
=\left(\frac{R/\PP_{i_1}}{(\QQ'\bigcap B_{i_1)}/\PP_{i_1}}\right)_{S_1}
=\left(\frac{B_{i_1}}{\QQ'\bigcap B_{i_1}}\right)_{S_1}=$$
$$
\frac{B_{i_1}S_1^{-1}}{(\QQ'\bigcap B_{i_1})S_1^{_1}}=
\frac{B_{i_1}S_1^{-1}}{\QQ'}=
\frac{(B_{i_1}S_1^{-1})/\PP'_{i_1i_2}}{\QQ'/\PP'_{i_1i_2}}=
\frac{B_{i_1i_2}}{\QQ'/\PP'_{i_1i_2}}$$
The ideal $\QQ=\pi_1^{-1}(\QQ'\bigcap B_{i_1})$ is prime, $H$-stable and with
zero intersection with $S_1$.
According Proposition 2.9, this proves that all prime factors of $B_{i_1i_2}$
are pure quantum.
The generators of $B_{i_1i_2}$ satisfy required conditions and
Lemma 3.3 is valid for
$B_{i_1i_2}$.
We localize this algebra on 
$x'_{n-i_1-i_2-1}$, get new system of generators in $\pi_2\pi_1(R)$ and
continue the process. 

After $k$ steps we get 
the semiprime $C$-algebra $B_{i_1,\ldots,i_k}$ generated by the elements
( here $i=i_1+\ldots+i_k$):
$$\{x_{n-i_1}^{\pm 1},(x'_{n-i_1-i_2-1})^{\pm 1}\ldots
,(x^{(k-2)}_{n-i_1-\ldots-i_{k-1}-(k-2)})^{\pm 1},x^{(k-1)}_{n-i-(k-1)},
x^{(k-1)}_{n-i-k},
\ldots,x^{(k-1)}_j,\ldots,\}.$$
There exists the $H$-stable semiprime ideal $\PP_{i_1\cdots i_k}$ in
$B_{i_1\cdots i_k}$ and there exists the denominator set
$S_{k-1}$ finitely generated by 
$x_{n-i_1},\ldots,x^{(k-2)}_{n-i_1-\ldots-i_{k-1}-(k-2)}, 
N_{k-1}=S_{k-1}\cap C$
in factor algebra such that
$$\left(\frac{R}{\PP_{i_1\cdots i_k}}\right)_{S_{k-1}}=B_{i_1\cdots i_k}$$
We assume that algebra $B_{i_1\cdots i_k}$ satisfies the conditions of 
Lemma 3.3 and shall prove that next algebra $B_{i_1\cdots i_{k+1}}$
also satisfies this conditions (mainly the condition that prime factors are
pure quantum).

By Lemma 3.3, we consider the localization
$B_{i_1,\ldots,i_k}S_k^{-1}$ where
the denominator set $S_k$ is generated by $S_{k-1}$, 
$ x^{(k-1)}_{n-i-(k-1)}$ and finitely generated $N_k\subset C^*$.
Denote $C_k=(C)_{N_k}$.
By Lemma 3.3, we get new system of generators in $B_{i_1,\ldots,i_k}S_k^{-1}$:
$$
\{x_{n-i_1}^{\pm 1},\ldots,(x^{(k-1)}_{n-i-(k-1)})^{\pm 1},
x^{(k)}_{n-i-k}
\ldots,x^{(k)}_j,\ldots,\}.$$
with $n-i-k\le j$.
More precisely, $x^{(k)}_j\in\pi_k\cdots\pi_1(R)$.
Consider the set of $i_{k+1}$ generators 
$$L_{k+1}=\{x_{n-i-k}^{(k)},\ldots,x^{(k)}_{n-i-i_{k+1}-(k-1)}\}$$
Let $\JJ^{(k)}_{i_1\cdots i_{k+1}}$ be the minimal semiprime ideal
in $B_{i_1\cdots i_k}S_k^{-1}$ containing $L_{k+1}$.
We denote by $X_k$ the set of minimal prime ideals
in $B_{i_1\cdots i_k}S_k^{-1}$ over $\JJ^{(k)}_{i_1\cdots i_{k+1}}$.
As above we decompose
$X_k=X_k'\bigcup X_k''\bigcup X_k'''$
where
$$X_k'=\{\QQ^{(k)}\in X_k:
\QQ^{(k)}\bigcap C_k= 0, x^{(k)}_{n-i-i_{k+1}-k}\notin\QQ^{(k)}\},$$ 
$$
X_2''=\{\QQ^{(k)}\in X_k:
\QQ^{(k)} \bigcap C_k\ne 0, x^{(k)}_{n-i-i_{k+1}-k}\notin\QQ^{(k)}\},$$
$$X_2'''=\{\QQ^{(k)}\in X_k: x^{(k)}_{n-i-i_{k+1}-k}\in\QQ^{(k)}\}.$$
We denote
$$\PP^{(k)}_{i_1\cdots i_{k+1}}=\bigcap_{\QQ^{(k)}\in X_k'}\QQ^{(k)}$$ 
$$
B_{i_1\ldots i_{k+1}}=\frac{B_{i_1\ldots i_k}S_k^{-1}}
{\PP^{(k)}_{i_1\ldots i_{k+1}}}
$$
$$
\pi_{k+1}:B_{i_1\ldots i_k}\mapsto B_{i_1\ldots i_k}S_k^{-1}
\mapsto B_{i_1\ldots i_{k+1}}$$
Notice that for every prime(semiprime) ideal $\II$ in $B_{i_1\cdots i_{k+1}}$
the ideal $\pi_k^{-1}(\II)$ is prime(semiprime) in $B_{i_1\cdots i_k}$.
The  $\PP^{(k)}_{i_1\ldots i_{k+1}}$ is $H$-stable semiprime.
It contains $L_{k+1}$,
does not contain
$x^{(k)}_{n-i-i_{k+1}-k} $ and have zero intersection with $C$.
We save the notation $S_k$ for the image $\pi_{k+1}(S_k)$. 
We define the  prime $H$-stable ideal
$$\PP_{i_1\cdots i_{k+1}}=\pi_1^{-1}\cdots\pi_k^{-1}
(\PP^{(k)}_{i_1\ldots i_{k+1}}\bigcap B_{i_1\cdots i_k})=
\pi_1^{-1}\cdots\pi_{k+1}^{-1}(0)$$
Since $\PP_{i_1\cdots i_{k}}=\pi_1^{-1}\cdots\pi_k^{-1}(0)$,
we obtain
$\PP_{i_1\cdots i_{k+1}}\supset\PP_{i_1\cdots i_{k}}$.
Note
$$
\pi_k\cdots\pi_1(\PP_{i_1\cdots i_{k+1}})S_{k-1}^{-1}=
\PP^{(k)}_{i_1\ldots i_{k+1}}\bigcap B_{i_1\cdots i_k}
$$
This implies
$$
\left(\frac{\PP_{i_1\cdots i_{k+1}}}{
\PP_{i_1\cdots i_k}}\right)_{S_{k-1}}=
\PP^{(k)}_{i_1\ldots i_{k+1}}\bigcap B_{i_1\cdots i_k}
$$
The above  ideal has zero intersection with subset $S_k$ in
$R/\PP_{i_1\cdots i_k}$.
We have
$$
(R/\PP_{i_1\cdots i_{k+1}})_{S_k}=
\frac{R/\PP_{i_1\cdots i_{k}}S_k^{-1}}{
(\PP_{i_1\cdots i_{k+1}}/\PP_{i_1\cdots i_{k}})S_k^{-1}}=
\frac{(R/\PP_{i_1\cdots i_{k}})S_{k-1}^{-1}S_k^{-1}}{
(\PP_{i_1\cdots i_{k+1}}/\PP_{i_1\cdots i_{k}})S_{k-1}^{-1}S_k^{-1}}=$$

$$\frac{B_{i_1\cdots i_k}S_k^{-1}}{
(\PP^{(k)}_{i_1\ldots i_{k+1}}\bigcap B_{i_1\cdots i_k})S_k^{-1}}=
\frac{B_{i_1\cdots i_k}S_k^{-1}}{\PP^{(k)}_{i_1\ldots i_{k+1}}}
= B_{i_1\cdots i_{k+1}}.
$$
If $\QQ^{(k)}\in X_k'$ then $\QQ^{(k)}$ is minimal prime ideal over 
$\PP^{(k)}_{i_1\ldots i_{k+1}}$ and
$\QQ^{(k)}/\PP^{(k)}_{i_1\ldots i_{k+1}}$ is minimal prime in
$B_{i_1\cdots i_{k+1}}$.
It follows that
the ideal 
$\QQ=\pi_1^{-1}\cdots\pi_k^{-1}
(\QQ^{(k)}\bigcap B_{i_1\cdots i_k})$
is minimal over $\PP_{i_1\cdots i_{k+1}}$ in $R$.
The ideal $\QQ$ has zero intersection with $S_k$ in factor algebra of $R$
over $\PP_{i_1\cdots i_k}$.
We also have
$$
\left(\frac{\QQ}{\PP_{i_1\cdots i_{k}}}\right)_{S_{k-1}}=
\QQ^{(k)}\bigcap B_{i_1\cdots i_k},
$$
$$
\left(\frac{\QQ}{\PP_{i_1\cdots i_{k+1}}}\right)_{S_{k-1}}=
\frac{\QQ^{(k)}\bigcap B_{i_1\cdots i_k}}
{\PP^{(k)}_{i_1\cdots i_{k+1}}\bigcap B_{i_1\cdots i_k}}
$$
$$
\left(\frac{\QQ}{\PP_{i_1\cdots i_{k+1}}}\right)_{S_{k}}=
\frac{(\QQ^{(k)}\bigcap B_{i_1\cdots i_k})S_k^{-1}}
{(\PP^{(k)}_{i_1\cdots i_{k+1}}\bigcap B_{i_1\cdots i_k})S_k^{-1}}=
\QQ^{(k)}/\PP^{(k)}_{i_1\cdots i_{k+1}}
$$
Now we are able to calculate prime factors
$$
(R/Q)_{S_k}= 
\frac{(R/\PP_{i_1\cdots i_{k+1}})S_k^{-1}}{(Q/\PP_{i_1\cdots i_{k+1}})S_k^{-1}}
=\frac{(B_{i_1\cdots i_k}S_k^{-1})/(\PP^{(k)}_{i_1\cdots i_{k+1}}) }
{\QQ^{(k)}/\PP^{(k)}_{i_1\cdots i_{k+1}}}=
\frac{B_{i_1\cdots i_{k+1}}}{\QQ^{(k)}/\PP^{(k)}_{i_1\cdots i_{k+1}}}
$$
Similar to case $k=2$ algebra $R/Q$ is pure quantum.
Then all prime factors of $B_{i_1\ldots i_{k+1}}$
over minimal ideals are pure quantum.
This proves that we may apply Lemma 3.3 and localize this algebra 
 over
$x^{(k)}_{n-i-i_{k+1}-k}$ and continue the process.

The process completes if $n-i-i_{k+1}-k=0$.
If this equality is true, then we get the $C$-algebra $B_{i_1\ldots i_{k+1}}$
generated by $q$-commuting elements
$$\{x_{n-i_1}^{\pm 1},\ldots,(x^{(k-1)}_{n-i-(k-1)})^{\pm 1}\}.$$
We obtain that the algebra
$$(R/\PP_{i_1\cdots i_{k+1}})_{S_k}= B_{i_1\cdots i_{k+1}}$$ 
is a factor algebra
of a algebra on twisted Laurent polynomials.
This completes Step 1.\\
{\bf Step 2}. Let $\II$ be a prime ideal with zero intersection
with $C$. 
Suppose that $\II\supset L_1$ for some $i_1$ and $x_{n-i_1}\notin \II$.
Then $\II$ contains one of the prime ideals
of the set  $X_1$.
Indeed $\II$ contains a unique
prime ideal $Q\in X_1'$ .
Then $\II\supset\PP_{i_1}$.
We cut $\PP_{i_1}$ and consider the 
ideal
$$\II'=\II/\PP_{i_1}$$ as a prime ideal in $R/\PP_{i_1}$.
The ideal $\II'$ has zero intersection with $S_1$.
We consider localization $S_1$ of $B_{i_1}=R/\PP_{i_1}$
over $x_{n-i_1}$. We get the ideal $\II'S_1^{-1}$
in $B_{i_1}S_1^{-1}$.
If $\II'S_1^{-1}$ 
contains some $L_2$,
do not contain $x'_{n-i_1-i_2-1}$ then $\II'S_1^{-1}\supset\PP'_{i_1i_2}$.
We continue.
After $k$-steps we get
$\II^{(k)}S_k^{-1}\supset \PP^{(k)}_{i_1\ldots i_{k+1}}$.
and $\II^{(k)}\bigcap S_k=\emptyset$.
We have 
$$
\II= \pi_1^{-1}\pi_2^{-1}\cdots\pi_k^{-1}
(\II^{(k)}S_k^{-1}\bigcap B_{i_1\cdots i_k})
\supset  \pi_1^{-1}\pi_2^{-1}\cdots\pi_{k}^{-1}
(\PP^{(k)}_{i_1\ldots i_{k+1}}\bigcap B_{i_1\cdots i_k})=
\PP_{i_1\ldots i_{k+1}} 
$$
and 
$$
S_k\bigcap \II\mod(\PP_{i_1\ldots i_{k+1}} )=
\emptyset.
$$
$\Box$\\
{\bf Remark}.
The number of elements of $\M$ is less o equal to $2^n$
The equality holds for the case $R$ is an algebra 
$C_Q[x_1, x_2,\ldots, x_n, k_1^{\pm},\ldots,k_m^{\pm}]$
of twisted polynomials.

Denote $\Spec_0(R)=\{\II\in\Spec(R), \II\bigcap C= 0\}$.
We equip $\Spec_0(R)$ with induced Zariski topology.\\
{\bf Corollary 3.5}.
There exists finite stratification of $\Spec_0(R)$ by the system
$\M$ of locally closed subsets in Jacobson topology.\\
{\bf Proof}. Consider the subset
$$
\M_{i_1\ldots i_{k+1}}=\{\II\in\Spec_0(R): 
\II\supset \PP_{i_1\ldots i_{k+1}}\quad{\mbox{and}}\quad
S_k\bigcap\II\mod(\PP_{i_1\ldots i_{k+1}})=\emptyset\}
$$
The set $S_k\in R/\PP_{i_1\ldots i_{k+1}}$ generated by $k$ elements, say
$y_1,\ldots,y_k$. Denote by $y'_1,\ldots,y'_k$ their 
representatives in inverse image.
The set $\M_{i_1\ldots i_{k+1}}$ is an intersection of Jacobson closed set
$$
\{\II\in\Spec_0(R): 
\II\supset \PP_{i_1\ldots i_{k+1}}\}$$
with the intersection of Jacobson open sets 
$
\{\II\in\Spec_0(R): y'_i\notin\II\}$.
This proves that $\M_{i_1\ldots i_{k+1}}$ is locally closed.
The rest follows from the Proposition.\\
{\bf Corollary 3.6}. Let $\II$ be a prime ideal in quantum solvable 
algebra $R$ obeying Q1-Q4 and $\II\bigcap C=0$. Denote $B=R/\II$.
Then there exist a denominator set $S$ in $B$ generated by $q$-commuting elements,
such that the ring $B_S$ is isomorphic to the factor ring of the algebra
of twisted Laurent polynomials.\\
{\bf Proof.} The statement 2) of the Proposition 3.4
implies that the exists the homomorphism
$$
\phi:\left(\frac{R}{\PP_{i_1\ldots i_{k+1}}}\right)_{S_k}
\mapsto 
\left(\frac{R}{\II}\right)_{S_k}.$$
The first algebra is the algebra of twisted Laurent polynomials.
This completes the proof.$\Box$\\
{\bf Corollary 3.7}.
$\II$, $R$ as Corollary 3.6.
Then $R/\II$ is pure quantum.
\\
{\bf Proof}. The algebra of twisted Laurent polynomials is pure quantum
[AD].\\ 
{\bf Corollary 3.8}. $\II$,$R$ as in Corollary 3.6.
The skew field of fractions of $R/\II$ is isomorphic to the skew field
of fractions of the algebra of twisted polynomials.\\
{\bf Proof}. 
We can restrict to the case  $\II$ 
is a ideal in the algebra of twisted Laurent polynomials
$B= C_Q[x_1^{\pm 1},\ldots, x_n^{\pm 1}]$.
Well known that 
every ideal in $B$ is generated by the intesection with its center.
The center $\Center(B)$ is generated by monomials
$x_1^{m_1}\cdots x_n^{m_n}$ where $\overline{m}=(m_1,\ldots,m_n)$
is an integer solution of system of equations
$$q_{11}^{m_1}\cdots q_{in}^{m_n}=1,\quad i\in\overline{1,n}.$$
We recall that the group $\Gamma$ is torsionfree.
The abelian group $G=\{\overline{m}\}$ of solutions of (2.9) is the subroup
in $\Z^n$. There exists compatible basises.
There exists the knew set of $q$-commuting generators 
$Y_1,\ldots,Y_m,Y_{m+1},\ldots,Y_n$
in $B$ such that $\Center(B)$ is generated by $Y_{m+1},\ldots,Y_n$.
It follows that $B/\II$ is generated by $Y_1,\ldots,Y_m$
over its center.
The skew field of fractions is the skew field of fractions of algebra
of twisted polynomials.$\Box$\\
{\bf Corollary 3.9}. Let $R$ be a quantum solvable algebra
obeying Conditions Q1-Q4.
As in Section 1 we denote $F=Fract(C)$. Consider the $F$-algebra 
$R_F=R\otimes_CF$. The Proposition 3.4 and Corollaries 3.5-3.7
is true for $R_F$ after replacing $R$ by $R_F$, $\Spec_0(R)$ by $\Spec(R_F)$,
a prime $R$-ideal with zero intersection with $C$ by
a prime $R_F$-ideal.\\
{\bf Proof}. Consider the denominator set $T=C^*$ in $R$.
The maps $\II\mapsto \II T^{-1}=\JJ\in\Spec(R_F)$ and 
$\JJ\mapsto\JJ\cap R=\II\in\Spec_0(F)$
establishes one to one correspondense between
$\Spec_0(R)$ and $\Spec(R_F)$. The algebras $(R/\II)_T$
and $R_F/\JJ$ are isomorphic. This completes the proof.
$\Box$

Our next goal is to show how Teorem 3.4 and its Corollaries
work in specialzations of $C$.
We consider the chain of ideals 
$\JJ^{(m)}_{i_1\ldots i_{m+1}}$, $m\in\overline{0,k}$ in
$B_{i_1\ldots\i_m}S_m^{-1}$ (see the Proof of Theorem 3.4).
For each prime ideal $\QQ^{(m)}\in X''_{m+1}$ we consider its
nonzero (by definition) intersection $\QQ^{(k)}\bigcap C$. 
Denote by $Y$ the union of all these prime ideals
that appears in construction of all standart ideals
$\PP_{i_1\cdots\i_{k+1}}$. The set $\overline{Y}$ is 
Zariski-closure of $Y$.
Consider also the set $S_k$ and $N_k=S_k\bigcap C$ (see Lemma 3.3
and Proof of Theorem 3.4). The set $N_k$ is finitely generated
by some elements, say  $g_1,\ldots,g_s$.
Consider the nonempty Zariski-open set  
$$\OO=\{\LL\in\Spec(C): \LL\notin\overline{Y}, 
\LL\bigcap\{g_1,\ldots,g_s\}=\emptyset\}.$$
Let $\la$ be a specialization $\la: C\mapsto K$ with
$\Ker(\la)=\LL$ is a maximal ideal in $C$ and $K=C/\LL$ is a field.
By definition of $C$, the field $K$ has zero characteristic. 
Denote $R_\la=R/(R\LL)$ and $\pi_\la:R\mapsto R_\la$ specialization of
$R$. For any ideal $\II_\la$ in $R_\la$ we denote 
$\II=\pi_\la^{-1}(\II_\la)$
in $R$.\\
{\bf Theorem 3.10}. Suppose that $R$ is quantum solvable algebra
obeying Conditions Q1-Q4 and $\LL\notin\OO$ then\\
1) for every completely prime ideal $\II_\la$ of $R_\la$ there exists
a unique standart ideal
$\PP_{i_1\cdots i_{k+1}}$ such that $\II\supset\PP_{i_1,\ldots i_{k+1}}$ and 
$S_k\bigcap(\II\mod(\PP_{i_1,\ldots i_{k+1}}))= \emptyset$;\\
2) the subsets
$$
\M^\la_{i_1\ldots i_{k+1}}=\{\II_\la\in\Spec(R_\la): 
\II\supset \PP_{i_1\ldots i_{k+1}}\quad{\mbox{and}}\quad
S_k\bigcap\II\mod(\PP_{i_1\ldots i_{k+1}})=\emptyset\}
$$
 are locally-closed in Jacobson topology of $\Spec(R_\la)$ and form
a stratification of $\Spec(R_\la)$;\\
3) If $\II_\la\in\M^\la_{i_1\ldots i_{k+1}}$, then the ring
$(R_\la/\II_\la)_{S_k}$ is a factor ring of the algebra of twisted
Laurent polynomials;\\
4) Denote $\la(q_{ij})=\la_{ij}$. If group
generated by $\la_{ij}$ is torsionfree, then $\Fract(R_\la/\II_\la)$
is isomorphic to the skew field of fractions of the algebra of twisted
polynomials.\\
{\bf Proof}. Let $\II_\la$ be a completely prime ideal in $R_\la$
and $\II$ the corresponding ideal in $R$.
Suppose that $\II$ contains $L_1$ and does not contain
$x_{n-i_1}$ (see Proof of Theorem 3.4).
Then $\II\supset \JJ_{i_1}$. It follows that $\II\supset\QQ\in X_1'
\bigcup X_1''$. If $\QQ\in X_1'$, then 
$\LL=\II\bigcap C\supset \QQ\bigcap C$. The last ideal lie in $Y$ and
$\LL\notin \overline{Y}$. We get a contradiction. It follows
that $\II$ contains some $\QQ\in X_1'$. Whence $\II\supset \PP_{i_1}$.
Notice that $\LL=\II\bigcap C$ does not contain the generators
$g_1,\ldots,g_s$ of $N_k$. Since 
$N_k\supset N_{k-1}\supset\cdots \supset N_1$ and $\II$ is completely prime, 
the ideal $\II$ has empty intersection with 
$N_1$ and whence $\II'=\II/\PP_{i_1}$ has empty intersection with $S_1$.
We consider localization $\II'S_1^{-1}$ as an ideal in $B_{i_1}$.
The proof of 1) may be completed similar to Step 2 in Theorem 3.4.
The statements 2),3,4) are proved as Corollaries 3.5-8.
$\Box$\\   
{\bf Corollary 3.11}. Suppose that $R$ be as in Theorem 3.10
and $C$ is the algebra of Laurent polynomials $\Q[q,q^{-1}]$
over the field of rational numbers.
Consider $R_\C= R\otimes_\Q\C$. 
Let $\la\in \C$ be a transcendental number.
Consider specialization $q\mapsto\la$ and $\Ker(\la)=\LL=\C[q,q^{-1}](q-\la)$.
We assert that $\LL\in\OO$ and whence all statements
of Theorem 3.10 is valid for specrum of specialzation $\Spec(R_\la)$.\\
{\bf Proof}. If $\LL$ not in $\OO$, then $\LL$ contains some
nonzero rational polynomial lying in  
$\QQ^{(k)}\bigcap C\in Y$ or $N_k$.
This polynomial becomes zero after reduction $\mod(q-\la)$.
This contradicts that $\la$ is a transcendental number.$\Box$.  \\
{\bf Example 3.12}.
Consider the quantum solvable algebra $R$ generated by 
$x,y$ over $C=\C[q,q^{-1}]$ with the single relation 
$xy-qyx=f(q)$. 
Here $f(q)$ is a complex polynomial with roots
$\alpha_1,\ldots,\alpha_n$.
Denote $u=xy-yx=(q-1)yx+f(q)$. There holds $uy=qyx$.
The stratification of $\Spec_0(R)$ (the prime ideals with zero
intersection with C) consists of 
$$\M_1= \{\II\in \Spec_0(R): y,u\notin \II\},$$ 
$$\M_2= \{\II\in \Spec_0(R): y\notin \II,u\in\II\}.$$
The ideal in $\Spec_0(R)$ can not contain $y$ (and $x$ also).
For any ideal $\II$ in $\Spec(R)-\Spec_0(R)$ we consider
$\II\bigcap C=C(q-\la)$.
Identifying here the prime ideals of $C$ with complex numbers,
we get  $\OO=\C-\{1,\alpha_1,\ldots,\alpha_n\}$. 
For $\la\in\OO$ we get statification of $\Spec(R_\la)$
as follows
$$\M_1^\la= \{\II_\la\in \Spec(R_\la): y,u\notin \II_\la\},$$ 
$$\M_2^\la= \{\II_\la\in \Spec(R_\la): y\notin \II_\la,u\in\II_\la\}.$$
But for $\la=\alpha_i\notin\OO$ we have the prime ideal 
$\II_i=<y, q-\alpha_i>$ which don't lie in $\M_1^\la$,$\M_2^\la$.     
In the case $\la=1\notin\OO$ and $c=f(1)\ne 0$ we get 
$xy-yx=c$. We get the Weyl algebra. The $\Fract(R)$
is not isomorphic to skew field of twisted rational functions.
We see that some statements of Theorem 3.10 may be false
if $\la\notin\OO$.

\end{document}